
\tolerance=10000
\raggedbottom

\baselineskip=15pt
\parskip=1\jot

\def\sk{\vskip 3\jot}

\def\heading#1{\vskip3\jot{\noindent\bf #1}}
\def\label#1{{\noindent\it #1}}
\def\QED{\hbox{\rlap{$\sqcap$}$\sqcup$}}


\def\ref#1;#2;#3;#4;#5.{\item{[#1]} #2,#3,{\it #4},#5.}
\def\refinbook#1;#2;#3;#4;#5;#6.{\item{[#1]} #2, #3, #4, {\it #5},#6.} 
\def\refbook#1;#2;#3;#4.{\item{[#1]} #2,{\it #3},#4.}


\def\({\bigl(}
\def\){\bigr)}


\def\al{\alpha}
\def\be{\beta}
\def\ga{\gamma}

\def\ze{\zeta}

\def\rh{\varrho}
\def\si{\sigma}
\def\ta{\tau}

\def\ch{\chi}

\def\Ga{\Gamma}



\def\abs#1{\bigg\vert#1\bigg\vert}
\def\norm#1{\|#1\|}

\def\lcm{{\rm lcm}}

\def\Not{N^{1/12}}
\def\SN{\{1, \ldots, N\}}

\def\Ex{{\rm Ex}}
\def\Var{{\rm Var}}

{
\pageno=0
\nopagenumbers
\rightline{\tt fleming.moments.tex}
\vskip1in

\centerline{\bf Large Deviations and Moments  for  the}
\centerline{\bf Euler Characteristic of a Random Surface}
\vskip0.5in

\centerline{Kevin Fleming}
\centerline{\tt kfleming@hmc.edu}
\vskip0.5in

\centerline{Nicholas Pippenger}
\centerline{\tt njp@hmc.edu}
\vskip0.5in

\centerline{Department of Mathematics}
\centerline{Harvey Mudd College}
\centerline{1250 Dartmouth Avenue}
\centerline{Claremont, CA 91711}
\vskip0.5in

\noindent{\bf Abstract:}
We study random surfaces constructed by glueing together $N/k$ filled $k$-gons along their edges,
with all $(N-1)!! = (N-1)(N-3)\cdots3\cdot 1$ pairings of the edges being equally likely.
(We assume that $\lcm\{2,k\}$ divides $N$.)
The Euler characteristic of the resulting surface is related to the number of cycles in a certain random permutation of $\SN$.
Gamburd has shown that when $2\,\lcm\{2,k\}$ divides $N$,
the distribution of this random permutation converges to that of the uniform distribution on the alternating group $A_N$ in the total-variation distance as $N\to\infty$.
We obtain large-deviations bounds for the number of cycles that, together with Gamburd's result,
allow us to derive sharp estimates for the moments of the number of cycles.
These estimates allow us to confirm certain cases of conjectures made by Pippenger and Schleich.
\vfill\eject
}

\heading{1.  Introduction}

The random surfaces that we deal with are defined as follows.
Fix $k\ge 3$ and let $N$ be divisible by $\lcm\{2,k\}$.
Take $N/k$ filled $k$-gons, and identify the $N$ edges of their boundaries in pairs to 
obtain an orientable surface, with all $(N-1)!! = (N-1)(N-3)\cdots3\cdot 1$ pairings  being equally likely.
We are interested in the Euler characteristic $\ch$ of the resulting surface.
Since the surface pattern has $N/k$ faces and $N/2$ edges (after identification),
we have $\ch = V - N/2 + N/k$, where $V$ is the number of vertices (after identification) in the surface pattern.
We shall thus focus our attention on the random variable $V$.

These random surfaces were first studied, for $k=3$, by Pippenger and Schleich [P],
who showed that
$$\eqalign{
\Ex[V] &= \log N + O(1), \cr
\Var[V] &= O(\log N), \cr
}$$
and conjectured that
$$\eqalignno{
\Ex[V] &= \log N + \ga  + o(1), &(1.1)\cr
\Var[V] &= \log N + \ga  - {\pi^2\over 6} + o(1), &(1.2)\cr
}$$
where $\ga = 0.5772\ldots$ is Euler's constant and $\pi = 3.14159\ldots$ is the circular ratio.

The random variable $V$ can be interpreted as the number of cycles in a random permutation.
Let $\al$ be a random permutation of $\SN$ uniformly distributed on the conjugacy class
$[2^{N/2}]$ (the class of permutations all of whose cycles have length $2$).
Let $\be$ be a fixed permutation in the conjugacy class $[k^{N/k}]$ (the class of permutations
all of whose cycles have length $k$).
Then the distribution of $V$ is the same as that of the number of cycles in the random permutation
$\al\be$.
(The elements of $\SN$ correspond to the edges of the edges of the $k$-gons, the permutation $\al$ corresponds to the pairing of the edges, and the permutation $\be$ corresponds to the cyclic ordering of the edges around each $k$-gon.)

Let us focus for now on the distribution of the cycle structure of the random permutation $\al\be$.
We denote by $\al\be$ the permutation obtained by first performing $\al$, then performing $\be$, but the
opposite convention (first $\be$, then $\al$) would yield the same distribution of the cycle structure: since $\al\be$ is conjugate to $\be\al = \be(\al\be)\be^{-1}$, they have the same cycle structure.
Gamburd [G], in the theorem cited below, assumes that $\be$ is uniformly distributed on the conjugacy class $[k^{N/k}]$, rather than being a fixed permutation from this class, but this change also yields the same distribution of the cycle structure: if $\be$ is a fixed permutation in $[k^{N/k}]$, and $\pi$ is uniformly distributed in the symmetric group $S_N$, then $\pi\be\pi^{-1}$ is uniformly distributed on $[k^{N/k}]$, but $\al(\pi\be\pi^{-1})$ is conjugate to 
$\pi^{-1}\(\al(\pi\be\pi^{-1})\)\pi = (\pi^{-1}\al\pi)\be$, which has the same distribution  as
$\al\be$, because $\pi^{-1}\al\pi$ has the same distribution as $\al$.
This argument also shows that the choice of the fixed permutation $\be$ in $[k^{N/k}]$ does not affect the
distribution of the cycle structure.
Finally, we observe that the distribution of the cycle structure of $\al\be$ actually determines the distribution of $\al\be$, since the distribution of $\al\be$ is constant on conjugacy classes:
If $\rh$ and $\pi\rh\pi^{-1}$ are conjugate permutations, then 
$\Pr[\al\be = \pi\rh\pi^{-1}] = \Pr[\pi^{-1}(\al\be)\pi = \rh] 
= \Pr[(\pi^{-1}\al\pi)(\pi^{-1}\al\pi)] = \Pr[\al\be = \rh]$, since $\pi^{-1}\al\pi$ has the same distribution as 
$\al$ and $\pi^{-1}\be\pi$ belongs to the same conjugacy class, $[k^{N/k}]$, as $\be$.

If $\pi$ is a random variable distributed on the symmetric group $S_N$, we shall denote by
$P_\pi$ its probability distribution, the function that maps each event $E\subseteq S_N$
to its probability $P_\pi[E] = \Pr[\pi\in E]$.
If $\rh$ is another such random variable, we shall denote by $\norm{P_\pi - P_\rh}$ the
{\it total-variation distance\/} between $\pi$ and $\rh$, that is, the maximum, over all events 
$E\subseteq S_N$, of the  absolute value of the
difference between probabilities assigned to $E$ by $P_\pi$ and $P_\rh$:
$$\norm{P_\pi - P_\rh} = \max_{E\subseteq S_N} \abs{P_\pi[E] - P_\rh[E]}.$$

\label{Theorem 1.1:} (Gamburd [G])
Fix $k\ge 3$, and let $N$ tend to infinity through values divisible by $2\,\lcm\{2,k\}$.
Let $\al$ be a random permutation uniformly distributed on the conjugacy class $[2^{N/2}]$ of 
the symmetric group $S_N$, let $\be$ be a fixed permutation in $[k^{N/k}]$,
and let $\ta$ be uniformly distributed on the alternating group $A_N$.
Then the distribution of $\al\be$ tends to that of $\ta$ in total-variation distance as $N\to\infty$.
More specifically,
$$\norm{P_{\al\be} - P_\ta} = O\left({1\over \Not}\right).$$
(Gamburd assumes that $\be$ is uniformly distributed in $[k^{N/k}]$, but as we have noted above,
this does not affect the distribution of $\al\be$.
Gamburd also assumes only that $N$ is divisible by $k$; but $N$ must be even
for $\al$ to be chosen from $[2^{N/2}]$, and
unless $N$ is also divisible by 
$2k$ when $k$ is even, and by $4k$ when $k$ is odd,
$\al\be$ will be distributed not on $A_N$ but on its coset in $S_N$ comprising the odd permutations,
hence our requirement that $N$ be divisible by $2\,\lcm\{2,k\}$.
In his proof of this theorem, Gamburd derives the estimate $\norm{P_{\al\be} - P_\ta}^2 = O(1/N^{1/6})$,
from which we obtain our estimate by taking square-roots on both sides.)

We are interested in the random variable $V$ of vertices in the surface pattern of a random surface,
which has the same distribution as the number $C_{\al\be}$ of cycles in the random permutation 
$\al\be$.
The conclusion of Theorem 1.1 is not very useful for estimating the moments of 
$C_{\al\be}$.
If, for example, $\pi$ is the identity in $A_N$ with probability $1/\Not$ and equal to $\ta$ with the remaining probability $1 - 1/\Not$, then $\norm{P_\pi - P_\ta} = O(1/\Not)$, but 
$\Ex[C_\pi] = N^{11/12} + O(\log N)$, while $\Ex[C_\ta] = \log N + O(1)$.
This example, however, suggests that what is needed is a large-deviations bound for $C_{\al\be}$,
and such a bound is provided by the following theorem.

\label{Theorem 1.2:}
Fix $k\ge 3$, and let $N$ tend to infinity through values divisible by $\lcm\{2,k\}$.
Let $\al$ be a random permutation uniformly distributed on the conjugacy class $[2^{N/2}]$ of 
the symmetric group $S_N$ and let $\be$ be a fixed permutation in $[k^{N/k}]$.
Then
$$\Pr[C_{\al\be} \ge t] = O\left(\left({2\over 3}\right)^{t/2}\,N\right).$$
(For this theorem we need only the weaker assumption that $N$ is divisible by $\lcm\{2,k\}$, rather than
by $2\,\lcm\{2,k\}$.)

Theorems 1.1 and 1.2 together allow us to relate the moments of $C_{\al\be}$ to those of $C_\ta$,
the number of cycles in a permutation $\ta$ uniformly distributed on the alternating group $A_N$.

\label{Theorem 1.3:}
Fix $k\ge 3$, and let $N$ tend to infinity through values divisible by $2\,\lcm\{2,k\}$.
Let $\al$ be a random permutation uniformly distributed on the conjugacy class $[2^{N/2}]$ of 
the symmetric group $S_N$, let $\be$ be a fixed permutation in $[k^{N/k}]$,
and let $\ta$ be uniformly distributed on the alternating group $A_N$.
Let $p(x)$ be a  polynomial of degree $l$ that is non-negative and non-decreasing for 
$x\ge 0$.
Then 
$$\Ex[p(C_{\al\be})]  = \Ex[p(C_\ta)] + O\left({(\log N)^l \over \Not}\right).$$
(A slightly weaker version of this theorem appears in the thesis of Fleming [F].
Since the proof of this theorem uses Theorem 1.1, we require that $N$ be divisible by 
$2\,\lcm\{2,k\}$, and not merely by $\lcm\{2,k\}$.)

Theorem 1.3, together with straightforward estimates of the moments of $C_\ta$, yields the following corollary.

\label{Corollary 1.4:}
Fix $k\ge 3$, and let $N$ tend to infinity through values divisible by $2\,\lcm\{2,k\}$.
Let $\al$ be a random permutation uniformly distributed on the conjugacy class $[2^{N/2}]$ of 
the symmetric group $S_N$ and let $\be$ be a fixed permutation in $[k^{N/k}]$.
Then
$$\eqalignno{
\Ex[C_{\al\be}] &= \log N + \ga + O\left({\log N \over \Not}\right) &(1.3) \cr
\Ex\big[(C_{\al\be} - \Ex[C_{\al\be}])^2\big]
&= \log N + \ga - {\pi^2\over 6} + O\left({(\log N)^2 \over \Not}\right) \cr
\Ex\big[(C_{\al\be} - \Ex[C_{\al\be}])^3\big]
&= \log N + \ga - {\pi^2\over 2} + 2\ze(3) + O\left({(\log N)^3 \over \Not}\right) \cr
}$$
and
$$\eqalign{  
\Ex\big[(C_{\al\be} - \Ex[C_{\al\be}])^4\big]
&= 3(\log N)^2 + (1 + 6\ga - \pi^2)\log N \cr
&\qquad + \left(\ga + 3\ga^3 - {7\pi^2\over 6} - \ga\pi^2 + 12\ze(3) + {\pi^4\over 60}\right) 
+ O\left({(\log N)^4 \over \Not}\right). \cr
}$$
This corollary confirms the conjectures (1.1) and (1.2)
made by Pippenger and Schleich [P] for $k=3$, at least in the case that $N$ is divisible by 
$2\,\lcm\{2,k\} = 12$.
\sk

\heading{2. Large Deviations}

In this section, we prove Theorem 1.2.
We begin by deriving an analog of Theorem 1.2 for the number of cycles $C_\si$ of a random permutation $\si$ uniformly distributed on the symmetric group $S_N$, as this will present the essential steps of our proof in a simpler context.

\label{Theorem 2.1:}
Let the random permutation $\si$ be uniformly distributed on the symmetric group $S_N$.
Then 
$$\Pr[C_\si \ge t] = O\left({N\over 2^t}\right).$$

\label{Proof:}
We begin by deriving the generating function $G_\si(x)$ for $C_\si$:
$$G_\si(x) = {x(x+1)\cdots(x+N-1)\over N!}. \eqno(2.1)$$
Although this identity is well known (the numerator on the right-hand side is the generating function
for the Stirling numbers of the first kind, which count permutations by their number of cycles;
see Stanley [S1, Proposition 1.3.4, p.~19] for three proofs),
we shall give an argument that can be elaborated to prove Theorem 1.2.

We shall describe a process for constructing the graph
$$\Ga_\si = \big\{\(i, \si(i)\) : i\in\SN\big\}$$
of $\si$ in $N$ steps.
We shall define a sequence $\Ga_{\si,m}$ of graphs ($0\le m\le N$) by setting
$\Ga_{\si,0} = \emptyset$, and for $0\le m\le N-1$ adding one pair to $\Ga_{\si,m}$ to
obtain $\Ga_{\si,m+1}$, ending with $\Ga_{\si,N} = \Ga_\si$.
For any $0\le m\le N$, each vertex in the graph $\Ga_{\si,m}$ will have in-degree at most $1$ and
out-degree at most $1$, and thus $\Ga_{\si,m}$ will be a disjoint union of paths and cycles.
In $\Ga_{\si,m}$ there will be $N-m$ paths, so in $\Ga_{\si,0}$ each vertex constitutes a path of length 
$0$, and in $\Ga_{\si,N}$ all vertices are in cycles.
Each path has a unique {\it head\/} (vertex with out-degree $0$) and a unique {\it tail\/}
(vertex with in-degree $0$).

At step $m$ ($0\le m\le N-1$), we choose a head $i$ arbitrarily from among the $N-m$ heads in 
$\Ga_{\si,m}$, then choose a tail $j$ uniformly from among the $N-m$ tails in $\Ga_{\si,m}$,
and finally set $\Ga_{\si,m+1} = \Ga_{\si,m} \cup \{(i,j)\}$.
This step closes a cycle if $j$ is the tail of the same path as $i$ is the head of.
Thus the probability that step $m$ closes a cycle is $1/(N-m)$.
It follows that $C_\si$ is the sum of $N$ independent random variables that assume the value $1$
with probabilities $1/N, 1/(N-1), \ldots, 1/2, 1$ and assume the value $0$ with the complementary 
probabilities $(N-1)/N, (N-2)/(N-1), \ldots, 1/2, 0$.
Since the generating function for a random variable that assumes the value $1$ with probability $q$ and
assumes the value $0$ with the complementary probability $1-q$ is $1-q + qx$, 
and since the generating function for a sum of independent random variables is the product of the generating functions of those random variables,
the generating function
for $C_\si$ is 
$$G_\si(x) = \((N-1)/N + x/N\)\((N-2)/(N-1) + x/(N-1)\)\cdots\(x/1\),$$
which is equivalent to (2.1).

From (2.1), we complete the proof of the large-deviations estimate in the usual way.
For $x\ge 1$ we have
$$\eqalign{
\Pr[C_\si \ge t]
& = \sum_{s\ge t} \Pr[C_\si = s] \cr
& \le {1\over x^t} \sum_{s\ge t} \Pr[C_\si = s] \, x^s \cr
& \le {G_\si(x) \over x^t}. \cr
}$$
Taking $x=2$, we obtain
$$\eqalign{
\Pr[C_\si \ge t]
&\le {G_\si(2) \over 2^t} \cr
&= {1\over 2^t} \, {2\cdot 3\cdots(N+1) \over N!} \cr
&= {N+1 \over 2^t}, \cr
}$$
which completes the proof of the theorem.
\QED

Later we shall need the analog of this theorem for the number $C_\ta$ of cycles of a random permutation uniformly distributed on the alternating group $A_N$.

\label{Theorem 2.2:}
Let the random permutation $\ta$ be uniformly distributed on the alternating group $A_N$,
where $N\ge 3$.
Then 
$$\Pr[C_\ta \ge t] = O\left({N\over 2^t}\right).$$

\label{Proof:}
We obtain the generating function $G_\ta$ for $C_\ta$ from (1.1) by bisection and renormalization:
$$\eqalignno{
G_\ta(x) 
&= {G_\si(x) + G_\si(-x)\over 2} \, {N!\over N!/2} \cr
&= { x(x+1)\cdots(x+N-1) + (-x)(-x+1)\cdots(-x+N-1)\over N!}. &(2.2) \cr
}$$
Proceeding as in the proof of Theorem 2.1, we obtain the same bound, since when we set
$x=2$, the third factor in the second term of the numerator of (2.2) vanishes, so 
$G_\ta(2) = G_\si(2)$.
\QED

\label{Proof of Theorem 1.2:}
We shall not be able to write an exact formula for the generating function $G_{\al\be}(x)$ for 
$C_{\al\be}$;
rather we shall derive the upper bound
$$G_{\al\be}(x) \le F(x^2), \eqno(2.3)$$
valid for $x\ge 1$, where
$$F(x) = {3x^2 (1+4x)(3+4x)\cdots(N-5+4x) \over (N-1)!!}, \eqno(2.4)$$
and where $M!! = M(M-2)\cdots 3\cdot 1$ for $M$ odd.

We shall describe a process for constructing the graph
$$\Ga_{\al\be} = \big\{\(i, \be\(\al(i)\)\) : i\in\SN\big\}$$
of $\al\be$ in $N/2$ steps.
We shall define a sequence $\Ga_{\al\be,m}$ of graphs ($0\le m\le N/2$) by setting
$\Ga_{\al\be,0} = \emptyset$, and for $0\le m\le N/2-1$ adding two pairs to $\Ga_{\al\be,m}$ to
obtain $\Ga_{\al\be,m+1}$, ending with $\Ga_{\al\be,N/2} = \Ga_\al\be$.

We fix $\be\in[k^{N/k}]$.
At step $m$ ($0\le m\le N/2-1$), we choose a head $i$ arbitrarily from among the $N-2m$ heads in 
$\Ga_{\al\be,m}$, then choose another head $j$ uniformly from among the $N-2m-1$ other heads in 
$\Ga_{\al\be,m}$,
and finally set $\Ga_{\al\be,m+1} = \Ga_{\si,m} \cup \{ \(i,\be(j)\), \(j,\be(i)\) \}$.

We now ask when cycles are created in this process.
A newly formed cycle can include just one of the two pairs added to the graph (we call this a 
{\it simple\/} closure) or both of the added pairs (we call this a {\it double\/} closure).
If vertex $h$ lies in a path, we shall denote by $H(h)$ the head of that path, and by $T(h)$ the tail of the path.
The functions $T$ and $H$ satisfy $H^2(h) = H(h)$ and  $T^2(h) = T(h)$, and also
$H\(T\(H(h)\)\) = H(h)$ and $T\(H\(T(h)\)\) = T(h)$.
It is easy to see that once $i$ has been chosen, there are just two values of $j$ that result in 
a simple closure: $j = \be^{-1}\(T(i)\)$ (which closes a cycle through the added pair $ \(i,\be(j)\)$)
and $j = H\(\be(i)\)$ (which closes a cycle through the added pair $ \(j,\be(i)\)$).

Double closures are more difficult to analyze.
We shall call a path a {\it quasi-cycle\/} if $T(h) = \be\(H(h)\)$ for any vertex $h$ on the path.
A double closure occurs when the two added pairs join two quasi-cycles into a single cycle.
The difficulty with focusing attention on this event is that there is no fixed bound on the number
of quasi-cycles that may be present at a given step, and thus no fixed bound on the number of values of 
$j$ that, together with a given choice of $i$, lead to a double closure.
For this reason, we shall transfer our attention to the creation of quasi-cycles, since the number of double closures during the process is bounded by the number of quasi-cycles created during the process (indeed, by half this number).
It is not hard to see that once $i$ has been chosen, there are at most two values of $j$ that result in the creation of a quasi-cycle: $j = \be^{-2}\(T(i)\)$ (when $\be^{-1}\(T(i)\)$ is a head) and
$j = H\(\be^2(i)\)$ (when $\be(i)$ is a head).

We shall say that step $m$ is {\it interesting\/} if either a simple closure occurs during that step or
a quasi-cycle is created during that step.
Since once $i$ has been chosen there are at most four values of $j$ (out of the $N-2m-1$ possible values) that will make the step interesting,
the number of interesting steps is stochastically dominated by the sum of $N/2$ independent random variables that assume the value $1$ with probabilities $4/(N-1), 4/(N-3), \ldots, 4/7, 4/5, 1, 1$
and assume the value $0$ with the complementary probabilities 
$(N-5)/(N-1), (N-7)/(N-3), \ldots, 3/7, 1/5, 0, 0$.
Thus the generating function for the number of interesting steps is at most
$$F(x) = \((N-5)/(N-1) + 4x/(N-1)\)\((N-7)/(N-3) + 4x/(N-3)\)\cdots
\(3/7 + 4x/7\)\(1/5 + 4/5\)\(3x/3\)\(x/1\)$$
for $x\ge 1$,
which is equivalent to (2.4).
Since the total number of cycles is at most twice the number of interesting steps, the generating
function $G_{\al\be}(x)$ for $C_{\al\be}$ is at most $F(x^2)$ for $x\ge 1$, which is equivalent to (2.3).

From (2.3) and (2.4), we complete the proof of the large-deviations estimate in the usual way.
For $x\ge 1$ we have
$$\eqalign{
\Pr[C_{\al\be} \ge t]
& = \sum_{s\ge t} \Pr[C_{\al\be} = s] \cr
& \le {1\over x^t} \sum_{s\ge t} \Pr[C_{\al\be} = s] \, x^s \cr
& \le {G_{\al\be}(x) \over x^t} \cr
& \le {F(x^2) \over x^t}. \cr
}$$
Taking $x=(3/2)^{1/2}$, we obtain
$$\eqalign{
\Pr[C_{\al\be} \ge t]
&\le \left({2\over 3}\right)^{t/2} \,F(3/2) \cr
&= {9\over 20} \,   \left({2\over 3}\right)^{t/2} \,{3\cdot 5\cdots(N+1) \over (N-1)!!} \cr
&\le \left({2\over 3}\right)^{t/2} \,(N+1), \cr
}$$
which completes the proof of the theorem.
\QED

We note that when $t / \log N\to\infty$, Theorems 1.2, 2.1 and 2.2 can all be improved by taking 
$x\approx t/\log N$.
But in our applications we shall have $t = O(\log N)$, and the versions we have given suffice
in this case.
\sk

\heading{3. Moments}

\label{Proof of Theorem 1.3:}
Summing by parts, we have
$$\Ex[C_{\al\be}] = p(0) + \sum_{1\le s\le N} \(p(s) - p(s-1)\)\,\Pr[C_{\al\be} \ge s].$$
Subtracting the analogous expression for $\Ex[C_\ta]$ and taking absolute values,
we obtain
$$\abs{\Ex[C_{\al\be}] - \Ex[C_\ta]}
\le  \sum_{1\le s\le N} \(p(s) - p(s-1)\)\,\abs{\Pr[C_{\al\be} \ge s] - \Pr[C_\ta \ge s]}$$
(where we have used the fact that $p(x)$ is non-decreasing for $x\ge 0$).
Breaking the sum at $s=t$ (where $t$ will be chosen later) yields
$$\eqalign{
\abs{\Ex[C_{\al\be}] - \Ex[C_\ta]}
&\le  \sum_{1\le s\le t} \(p(s) - p(s-1)\)\,\abs{\Pr[C_{\al\be} \ge s] - \Pr[C_\ta \ge s]} \cr
&\qquad +  \sum_{t <  s\le N} \(p(s) - p(s-1)\)\,\Pr[C_{\al\be} \ge s] \cr
&\qquad +  \sum_{t <  s\le N} \(p(s) - p(s-1)\)\,\Pr[C_\ta \ge s] \cr
&\le\norm{P_{\al\be} - P_\ta} \, p(t)
 +  \Pr[C_{\al\be} > t]\,p(N) + \Pr[C_\ta > t]\,p(N) \cr
}$$
(where we have used the fact that $p(x)$ is non-negative and non-decreasing
for $x\ge 0$).
Bounding the first term by Theorem 1.1, the second by Theorem 1.2 and the third by Theorem 2.2, we obtain
$$\abs{\Ex[C_{\al\be}] - \Ex[C_\ta]} = 
O\left({t^l \over \Not}\right) + O\left(\left({2\over 3}\right)^{t/2}\,N^{l+1}\right) 
+ O\left({N^{l+1}\over 2^t}\right)$$
(where we have used the fact that $p$ is of degree $l$).
Taking $t = \lceil 2(l+1)\log_{3/2} N\rceil$ completes the proof of the theorem.
\QED

To use Theorem 1.3, we must estimate the moments of $C_\ta$.
To  do that we shall first estimate the moments of $C_\si$.
The easiest moments to estimate are the ``factorial moments'',
$\Ex[C_\si^{\underline{l}}]$, where we have written $x^{\underline{l}}$ for the $l$-th  ``falling power''
$x(x-1)\cdots(x-l+1)$.

\label{Theorem 3.1:}
We have
$$\Ex[C_\si^{\underline{l}}] = (-1)^l \, Z\(-\ze_N(1), -\ze_N(2), \ldots, -\ze_N(l)\), \eqno(3.1)$$
where the polynomial $Z(g_1, g_2, \ldots, g_l)$ is defined by
$$Z_l(g_1, g_2, \ldots, g_l) = 
 \sum_{n_1 \ge 0, n_2 \ge 0, \ldots, n_l \ge 0 \atop
n_1 + 2n_2 + \cdots + ln_l = l} 
{l! \, g_1^{n_1} g_2^{n_2} \cdots g_l^{n_l} \over
n_1! 1^{n_1}\,n_2! 2^{n_2}\cdots n_l! l^{n_l}} \eqno(3.2)$$
is the ``cycle indicator'' for the symmetric group $S_l$,
and where $\ze_N(m)$ is defined by
$$\ze_N(m) = \sum_{1\le n\le N} {1\over n^m}.$$
(The sum $\ze_N(1)$, the $N$-th ``harmonic number'', is usually denoted $H_N$, but we have adopted a different notation for uniformity.)

\label{Proof:}
We have
$$\Ex[C_\si^{\underline{l}}] = {d^l G_\si(x)\over dx^l}\bigg\vert_{x=1}.$$
Each successive differentiation removes one of the factors $x, (x+1), \ldots, (x+N-1)$ from the numerator of $G_\si(x)$ (see (2.1)),
and $l$ differentiations remove $l$ distinct factors in some order,
giving rise to $N^{\underline{l}}$ terms, among which  ${N\choose l}$ distinct terms each appear $l!$ times.
When $x$ is set to $1$, the remaining factors in each term
cancel equal factors in the denominator, and we are then
left with the $l$ factors among $1, 2, \ldots, N$ in the denominator that correspond to those cancelled in the numerator.
Thus we have
$$\Ex[C_\si^{\underline{l}}] = l!\,e_l(1, 1/2, \ldots, 1/N),$$
where $e_l(x_1, x_2, \ldots, x_N)$ is the $l$-th ``elementary symmetric function'' of the variables
$x_,1 x_2, \ldots, x_N$:
$$e_l(x_1, x_2, \ldots, x_N) = \sum_{1\le i_1 < i_2 < \cdots  < i_l\le N} x_{i_1} x_{i_2} \cdots x_{i_l}.$$
We next express the $e_l$ in terms of  the ``power-sum symmetric functions'',
$$p_m(x_1, x_2, \ldots, x_N) = \sum_{1\le i\le N} x_i^m.$$
The relationship between the $e_l$ and the $p_m$ is given by
$$l!\,e_l = (-1)^l\,Z_l(-p_1, -p_2, \ldots, -p_l),$$
where the polynomial $Z_l(g_1, g_2, \ldots, g_l)$ is as given in (3.2)
(see Stanley [S2, Proposition 7.7.6, (7.23), p.~301 and Definition 7.24.1, p.~390],
Comtet [C, Chapter III, Exercise 9, pp.~158--159 and Equation [6e], p.~247], or
Riordan [R, Chapter 2, Exercise 27, p.~47 and Equation (2), p.~68].
Note that definitions of the cycle indicator differ among these authors
by a factor of the order of the group, in this case $l!$).
Since we have
$p_m(1,1/2,\ldots, 1/N) = \ze_N(m)$, we obtain (3.1).
\QED

We next show that the factorial moments of $C_\ta$ have the same asymptotic behavior as those of 
$C_\si$.

\label{Theorem 3.2:}
$$\Ex[C_\ta^{\underline{l}}] = \Ex[C_\si^{\underline{l}}] + O\left({(\log N)^{l-1} \over N}\right).$$

\label{Proof:}
From (2.2) we have $G_\ta(x) = G_\si(x) + G_\si(-x)$, from which we obtain
$$\Ex[C_\ta^{\underline{l}}] = \Ex[C_\si^{\underline{l}}] + 
{d^l G_\si(-x)\over dx^l}\bigg\vert_{x=1}.$$
Thus it will suffice to show
$${d^l G_\si(-x)\over dx^l}\bigg\vert_{x=1} = O\left({(\log N)^{l-1} \over N}\right). \eqno(3.3)$$
Substituting $x = -y$ and using the chain rule yields
$${d^l G_\si(-x)\over dx^l}\bigg\vert_{x=1} = (-1)^l\,{d^l G_\si(y)\over dy^l}\bigg\vert_{y=-1}. \eqno(3.4)$$
As in the proof of Theorem 3.1, the $l$ successive differentiations remove $l$ distinct factors from the
numerator of $G_\si(y)$.
But since $y$ is set to $-1$ after differentiation, the product of the remaining factors will vanish unless
the factor $y+1$ is removed by one of the $l$ differentiations.
Thus we have
$${d^l G_\si(y)\over dy^l}\bigg\vert_{y=-1} =
l\,{d^{l-1} \over dy^{l-1}}\,{y(y+2)\cdots(y+N-1)\over N!}\bigg\vert_{y=-1}. \eqno(3.5)$$
The product rule (for the product of $y$ and $(y+2)\cdots(y+N-1)$) now yields
$$\eqalignno{
{d^{l-1} \over dy^{l-1}}\,&{y(y+2)\cdots(y+N-1)\over N!}\bigg\vert_{y=-1} \cr
&=
\left((l-1)\,{d^{l-2} \over dy^{l-2}}\,{(y+2)\cdots(y+N-1)\over N!}
+ y\,{d^{l-1} \over dy^{l-1}}\,{(y+2)\cdots(y+N-1)\over N!}\right)\bigg\vert_{y=-1}. &(3.6)\cr
}$$
The remaining derivatives can now be expressed in terms of symmetric functions as in Theorem 3.1:
$${d^m \over dy^{m}}\,{(y+2)\cdots(y+N-1)\over N!}
= {m! \over N(N-1)} \, e_{m}\(1,1/2, \ldots, 1/(N-2)\). \eqno(3.7)$$
From
$$\ze_M(1) = \log M + \ga + O\left({1\over M}\right), \eqno(3.8)$$
where $\ga = 0.5772\ldots$ is Euler's constant, we obtain
$$m!\,e_m(1, 1/2, \ldots, 1/M) \le \ze_M(1)^m = O\((\log M)^m\).$$
Combining this estimate with (3.4) through (3.7) yields (3.3).
\QED

\label{Proof of Corollary 1.4:}
Since $Z_1(g_1) = g_1$, $Z_2(g_1,g_2) = g_2 + g_1^2$, 
$Z_3(g_1,g_2,g_3) = 2g_3 + 3g_1 g_2 + g_3^3$ and
$Z_4(g_1,g_2,g_3,g_4) = 6g_4 + 8g_3 g_1 + 3g_2^2 + 6g_2 g_1^2 + g_1^4$, Theorems 3.1 and 3.2
yield
$$\eqalign{
\Ex[C_\ta]
&= \ze_N(1) + O\left({1\over N}\right), \cr
\Ex[C_\ta^{\underline{2}}]
&= \ze_N(1)^2  - \ze_N(2) + O\left({\log N\over N}\right), \cr
\Ex[C_\ta^{\underline{3}}]
&= \ze_N(1)^3  - 3\ze_N(1)\,\ze_N(2) + 2\ze_N(3) + O\left({(\log N)^2 \over N}\right) \cr
}$$
and
$$\Ex[C_\ta^{\underline{4}}] = \ze_N(1)^4 - 6\ze_N(1)^2\,\ze_N(2)  + 8\ze_N(1)\,\ze_N(3)
+ 3\ze_N(2)^2 - 6\ze_N(4)$$
for the first four factorial moments of $C_\ta$.

From the factorial moments $\Ex[C_\ta^{\underline{l}}]$ of $C_\ta$, we can obtain the standard moments
$\Ex[C_\ta^{{l}}]$ via the expansion
$$\Ex[C_\ta^{{l}}] = \sum_{0\le m\le l} \left\{{l\atop m}\right\} \, \Ex[C_\ta^{\underline{m}}],$$
where the $\left\{{l\atop m}\right\}$ are the ``Stirling numbers of the second kind''
(see Stanley [S1, p.~209], or Comtet [C, Section 5.2, Theorem~B, p.~207])
and we can relate these moments to the standard moments $\Ex[C_{\al\be}^{{l}}]$ of $C_{\al\be}$ by 
Theorem 1.3.
This yields
$$\eqalignno{
\Ex[C_{\al\be}]
&= \ze_N(1) + O\left({\log N\over \Not}\right), &(3.9)\cr
\Ex[C_{\al\be}^{{2}}]
&= \ze_N(1)^2   + \ze_N(1) - \ze_N(2) + O\left({(\log N)^2\over \Not}\right), &(3.10) \cr
\Ex[C_{\al\be}^{{3}}]
&= \ze_N(1)^3  + 3\ze_N(1)^2 + \ze_N(1) - 3\ze_N(1)\,\ze_N(2) 
 - 3\ze_N(2) + 2\ze_N(3) + O\left({(\log N)^3 \over \Not}\right) &(3.11) \cr
}$$
and
$$\eqalignno{
\Ex[C_{\al\be}^4] 
&= \ze_1^4 + 6\ze_N(1)^3 + 7\ze_N(1)^2  - \ze_N(1)^2\,\ze_N(2) + \ze_N(1)
-18\ze_N(1)\,\ze_N(2) \cr
&\qquad + 8\ze_N(1)\,\ze_N(3)
+ 3\ze_N(2)^2 - 7\ze_N(2) + 12\ze_N(3) - 6\ze_N(4)
+ O\left({(\log N)^4 \over \Not}\right). &(3.12)\cr
}$$
From (3.8) and (3.9), we obtain (1.3) of Corollary 1.4.

From the standard moments of $C_{\al\be}$ we can obtain the central moments 
$\Ex\big[\(C_{\al\be} - \Ex[C_{\al\be}]\)^l\big]$ via the expansion
$$\Ex\big[\(C_{\al\be} - \Ex[C_{\al\be}]\)^l\big] = 
\sum_{0\le m\le l} {l\choose m} \Ex[C_{\al\be}^m] \, \(-\Ex[C_{\al\be}]\)^{l-m}.$$
Thus from (3.9) through (3.12) we obtain
$$\eqalign{
\Ex\big[(C_{\al\be} - \Ex[C_{\al\be}])^2\big]
&= \ze_N(1) - \ze_N(2) + O\left({(\log N)^2 \over \Not}\right), \cr
\Ex\big[(C_{\al\be} - \Ex[C_{\al\be}])^3\big]
&= \ze_N(1) - 3\ze_N(2) + 2\ze_N(3) + O\left({(\log N)^3 \over \Not}\right), \cr
}$$
and
$$\eqalign{
\Ex\big[(C_{\al\be} - \Ex[C_{\al\be}])^4\big]
&= 3\ze_N(1)^2 + \ze_N(1) - 6\ze_N(1)\,\ze_N(2)  \cr
&\qquad +  3\ze_N(2)^2 - 7\ze_N(2)
 + 12\ze_N(3) - 6\ze_N(4) + O\left({(\log N)^4 \over \Not}\right). \cr
}$$
Substituting (3.9) and
$$\ze_N(m) = \ze(m) + O\left({1\over N^{m-1}}\right)$$
for $m\ge 2$, where $\ze(m) = \sum_{n\ge 1} 1/n^m$ is Riemann's zeta function,
yields the remaining parts of Corollary 1.4, since $\ze(2) = \pi/6$ and $\ze(4) = \pi^4/90$
(see Brassoud [B, Appendix A.3, (A.37) and (A. 39), p.~285).
\sk

\heading{4. Acknowledgment}

The research reported here was supported
by Grant CCF 0430656 from the National Science Foundation.
\sk

\heading{5.  References}

\refbook B; D. M. Brassoud;
A Radical Approach to Real Analysis;
2nd Edition, Mathematical Association of America, 2007.

\refbook C; L. Comtet;
Advanced Combinatorics;
D.~Reidel Publishing Co., 1973.

\refbook F; K. Fleming;
Boundary Cycles in Random Triangulated Surfaces;
B.~S. Thesis, 2008, Department of Mathematics,
Harvey Mudd College, Claremont, CA 91711.

\ref G; A. Gamburd;
``Poisson-Dirichlet Distribution for Random Belyi Surfaces'';
Annals of Probability; 34:5 (2006) 1827--1848.


\ref P; N. Pippenger and K. Schleich;
``Topological Characteristics of Random Triangulated Surfaces'';
Random Structures and Algorithms; 28 (2006) 247--288.

\refbook R; J. Riordan;
An Introduction to Combinatorial Analysis;
John Wiley \& Sons, Inc., 1958 (reprinted by Dover, 2002).

\refbook S1; R. P. Stanley;
Enumerative Combinatorics;
v.~1, Cambridge University Press, 1997.

\refbook S2; R. P. Stanley;
Enumerative Combinatorics;
v.~2, Cambridge University Press, 1999.

\bye